\documentclass[a4paper,12pt,reqno]{article}
\usepackage[english]{babel}
\usepackage{amsmath,amsfonts,amssymb,amsthm,bbm}
\usepackage{graphics,epsfig,psfrag} 
\usepackage{paralist,subfigure,multirow,float}
\usepackage{hyperref} 
\usepackage[active]{srcltx} 
\usepackage[dvipsnames]{xcolor}
\usepackage{soul} 
\usepackage[latin1]{inputenc}
\usepackage[OT1]{fontenc}

\newtheorem{theorem}{Theorem}[section]

\newtheorem{definition}[theorem]{Definition}
\theoremstyle{definition}

 \usepackage{autonum} 

\DeclareMathOperator\inter{int}

\def\N{\mathbb{N}}
\def\Z{\mathbb{Z}}

\def\R{\mathbb{R}}

\let\O=\Omega
\let\e=\varepsilon
\def\epsilon{\varepsilon}

\let\t=\tilde
\let\ol=\overline

\let\.=\cdot
\let\0=\emptyset

\let\mc=\mathcal

\def\1{\mathbbm{1}}
\def\Sph{{\mathbb{S}}^{N-1}}
\newcommand{\su}[2]{\genfrac{}{}{0pt}{}{#1}{#2}}

\def\hat{\widehat}
\def\tilde{\widetilde}

\def\thm#1{Theorem~\ref{thm:#1}}
\def\seq#1{(#1_n)_{n\in\N}}

\newenvironment{formula}[1]{\begin{equation}\label{#1}}{\end{equation}\noindent}

\numberwithin{equation}{section}
\newcommand{\be}{\begin{equation}}
\newcommand{\ee}{\end{equation}}
\newcommand{\baa}{\begin{array}}
\newcommand{\eaa}{\end{array}}
\newcommand{\ba}{\begin{eqnarray}}
\newcommand{\ea}{\end{eqnarray}}
\def\Fi#1{\begin{formula}{#1}}
\def\Ff{\end{formula}\noindent}
\def\BS{\color{Bittersweet}}

\def\pe{principal eigenvalue}

\def\PTF{pulsating traveling front}
\def\ass{asymptotic set of spreading}
\def\W{\mc{W}}

\def\tinfty{\quad\text{ as }\;t\to+\infty}

\begin{document}
\date{}
\title{\bf{Freidlin--G\"artner formula and asymptotic profile in reaction-diffusion equations }}
\author{Luca Rossi$^{\hbox{\small{ c}}}$
	\vspace{7pt}\\
\footnotesize{$^{\hbox{c }}$SAPIENZA Univ Roma, Istituto ``G.~Castelnuovo'', Roma, Italy}}
\maketitle
	
\begin{abstract}
\noindent{
		We address the question of the large-time behavior of solutions to reaction-diffusion equations in periodic media. We start with the description of the asymptotic shape of the invasion set, which is characterized by the Freidlin--G\"artner formula. We outline a proof of the formula that holds true for general types of reaction terms. We then present some recent results, obtained in collaboration with H.~Guo and F.~Hamel, 
		for (weakly) bistable equations. 
		They include a regular version of the Freidlin--G\"artner formula and the convergence in profile towards pulsating traveling fronts 
		for solutions with either bounded or unbounded initial support.}
\end{abstract}
	


\section{Introduction}

Reaction-diffusion equations provide a fundamental framework for the study of
propagation phenomena in various fields of physical and life sciences. A central question is to understand
the large-time behavior of solutions arising from localized initial data.

In homogeneous environments, propagation occurs with a constant speed and the
asymptotic shape of the invaded region is a ball. In periodic media, however,
propagation becomes anisotropic: the speed depends on the direction, leading to
a family of spreading speeds and to a nontrivial asymptotic shape.

Another important source of complexity arises when considering initial data with
non-compact support. In this case, the geometry of the initial support shapes the \ass, 
and convergence to traveling fronts is more subtle.

In this note we will review the following topics:
\begin{itemize}
	\item expression of the \ass\ through the Freidlin--G\"artner formula;
	\vspace{-5pt}
	\item a unified proof of the formula taken from \cite{RossiFG} that applies to general reaction terms;
	\vspace{-5pt}
	\item some recent results from \cite{GHR1} about a ``regular'' Freidlin--G\"artner formula, as well as
	the convergence in profile towards pulsating traveling fronts.
\end{itemize}


\section{Benchmark: the homogeneous case}

Before turning to heterogeneous media, it is instructive to recall the classical
results in the homogeneous setting. These provide the intuition and the guiding
principles for the general case.

Consider the equation
\begin{equation}\label{eq:homogeneous}
	\partial_t u - \Delta u = f(u), \qquad t>0,\ x\in\R^N,
\end{equation}
where $f$ is a smooth function vanishing at $0$ and $1$, 
which represent, respectively, extinction and full invasion states.
We further assume that $f$ is of one of the following~types:
%
\begin{align}
\text{\em Monostable }\quad & f>0 \text{ in } (0,1)\\
& \text{a sub-class of which is the Fisher-KPP non-linearity}\\
&\text{satisfying in addition $f(s)\leq f'(0)s$ for $s\geq0$}\\[10pt]
\text{\em Ignition } \quad &
\exists \theta\in(0,1),\ \
\text{$f=0$ in $[0,\theta]$, $f>0 \text{ in } (\theta,1)$}\\[10pt]
\text{\em Bistable }\quad &
\exists\theta\in(0,1), \ \
\text{$f<0$ in $(0,\theta)$, $f>0 \text{ in } (\theta,1)$, with $\int_0^1 f>0$,}
\end{align}
and it is equal to $0$ outside $[0,1]$.
In population dynamics, the Ignition and Bistable cases refer to the weak and strong Allee effect respectively.

\subsection{Planar traveling fronts}

In all the cases: Monostable, Ignition and Bistable, the equation \eqref{eq:homogeneous} 
admits {\em planar traveling fronts}, that is, for a given direction $e\in\Sph$, solutions
of the form
\[
u(t,x)=\phi(x\.e-ct),
\]
with
\[
\phi(-\infty)=1, \qquad \phi(+\infty)=0.
\]
The quantity $c\in\R$ is called the {\em speed} of the front, 
and $\phi$ is called the {\em profile}. The latter satisfies the ODE $\phi''+c\phi'+f(\phi)=0$
in $\R$.
In the Monostable case, there exists a minimal speed $c^*>0$ such that
planar traveling fronts exist if and only if $c\geq c^*$.
In the Ignition and Bistable cases, there exists a unique possible speed $c=c^*>0$ for which a traveling fronts exists
(and is unique up to shift). In all the cases, we refer to $c^*$ as the {\em critical speed} of fronts.

\subsection{Spreading from compactly supported data}

A fundamental result due to Aronson and Weinberger describes the large-time behavior
of {\em invading} solutions emerging from compactly supported initial data, in all cases:
Monostable, Ignition and Bistable. Invasion means that the solution converges
locally uniformly to $1$ as $t\to+\infty$; see \cite{AronsonWeinberger} for 
some sufficient conditions.

\begin{theorem}[\cite{AronsonWeinberger}]\label{thm:homo}
	Let $u$ be an invading solution to \eqref{eq:homogeneous} with a compactly supported initial datum $u_0\geq0,\not\equiv0$.
	Then
	\[
	\lim_{t\to\infty}\Big( \inf_{|x|\le ct} u(t,x)\Big) = 1 \quad \text{if } c<c^*,
	\]
	and
	\[
	\lim_{t\to\infty}\Big( \sup_{|x|\ge ct} u(t,x)\Big) = 0 \quad \text{if } c>c^*.
	\]
\end{theorem}

	This result shows that propagation occurs in every direction with the same asymptotic speed $c^*$, which coincides with the minimal
	speed of planar fronts. More precisely, the upper level sets 
	\[
	E_\lambda(t)=\{x\in\R^N : u(t,x)>\lambda\},\qquad\lambda\in(0,1),
	\]
	of the solution expand as the ball of radius $c^* t$, up to a $o(t)$ term, thus, rescaling by $1/t$, one has
\[
\frac{1}{t} E_\lambda(t) \to B_{c^*}\tinfty,
\]
in the sense of the \emph{Hausdorff distance}, where $B_{c^*}$ is the ball of radius $c^*$.
The rescaled limit $B_{c^*}$ is referred to as the {\em \ass} (or invasion set) of the solution. 
This reveals that a symmetrization of the initial datum occurs along the evolution.
A result of Jones \cite{Jones} actually improves this conclusion, showing that the upper level set $E_\lambda(t)$ remains at bounded 
Hausdorff distance from some expanding ball (that a fortiori is of the form $B_{m(t)}$
with $m(t)=c^* t+o(t)$, by Theorem \ref{thm:homo}).
However, it is not true in general that such a distance converges to $0$ as $t\to+\infty$, see 
\cite{RossiJones,vsJones2,vsJones1}.

\subsection{Convergence to planar traveling fronts}\label{sec:profile-homo}

Heuristically, the reason why the solution with a compactly supported initial datum propagates with
the minimal speed $c^*$ of planar fronts, in any direction, is that it is ``becoming locally planar'' and therefore it is
approaching a planar traveling front. More precisely, 
for a solution $u$ as in Theorem \ref{thm:homo},
the following convergence holds for any $e\in\Sph$:
	\[
	u(t,x+m_e(t)e) \to \phi(x\.e)\tinfty\quad\text{ locally uniformly in $x\in\R^N$},
	\]
where $\phi$ is the profile of the minimal/unique (in the Monostable/Ignition or Bistable case) 
planar traveling front, 
and $m_e:\R^+\to\R$ is a suitable shift function.
We refer to~\cite{BHgtw,DM2,Ducrot,Gartner,RRR,Uchi-bi}
for the precise statements of this result for the different types of reaction term.
%
One knows from Theorem \ref{thm:homo} that $m_e(t)=c^*t+o(t)$.
It turns out that the amplitude of the 
``gap'' $c^*t-m_e(t)$ 
depends on the type of nonlinearity, as well as on the dimension~$N$, due to curvature effects.
There is a vast literature on this topic, especially in the Fisher-KPP case, where the problem can be 
tackled with both PDE and probabilistic techniques,
see \cite{Bram2,Ducrot,Gartner,HNRR1,RRR}.

\subsection{Limitations and motivation}

In homogeneous media:
\begin{itemize}
	\item propagation is isotropic,
	\item a single speed $c^*$ characterizes the dynamics,
	\item the asymptotic shape is a ball.
\end{itemize}

In periodic media, none of these properties remain true.


\section{Propagation in periodic heterogeneous media}\label{sec:FG}

Let us consider now reaction-diffusion equations in periodic heterogeneous media. 
These equations have the form
\begin{equation}\label{eq:main}
	\partial_t u - \nabla\cdot(A(x)\nabla u) - q(x)\cdot\nabla u = f(x,u),
	\qquad t>0,\ x\in\R^N.
\end{equation}
Throughout this note, we assume the following standing hypotheses.
The terms $A(x)$, $q(x)$ and $f(x,\cdot)$ are smooth and $\Z^N$-periodic in $x$, that is, 
$$\forall\,h\in\Z^N,\qquad A(\.+h)\equiv A,\qquad q(\.+h)\equiv q,\qquad f(\.+h,\.)\equiv f.$$
The matrix field $x\mapsto A(x)=(a_{ij}(x))_{1\le i,j \le N}$ is symmetric positive definite, 
the vector field $x\mapsto q(x)=(q_i(x))_{1\le i\le N}$ satisfies
\be\label{hypq}
\hbox{div}\, q=0 \hbox{ in $\R^N$},\qquad \int_{(0,1)^N}\!q(x)dx=0
\ee
and the function $f$ is smooth and satisfies
\be\label{f01}
f(x,0)=f(x,1)=0\ \ \hbox{ for all }\; x\in\R^N,
\ee
\be\label{fS}
\exists S\in(0,1),\ \text{ $f(x,\.)$ is nonincreasing in $[S,1]$}.
\ee
We complete equation \eqref{eq:main}
with an initial datum $0\leq u_0\leq 1$.


\subsection{The \ass}

In analogy with the homogeneous framework, our goal is to determine the {\em \ass} of
solutions of \eqref{eq:main}.

\begin{definition}\label{def:W}
We say that a solution $u$ to~\eqref{eq:main} admits an {\em \ass} $\W\subset\R^N$ if $\W$ coincides  with the interior of its closure and satisfies 
\Fi{ass-cpt}
\begin{cases}
\displaystyle\lim_{t\to+\infty}\,\Big(\min_{x\in C}u(t, tx)\Big)=1 & \text{for any compact set }C\subset\mc{W},\vspace{3pt}\\
\displaystyle\lim_{t\to+\infty}\,\Big(\max_{x\in C}u(t,tx)\Big)=0 & \text{for any closed set }C\subset\R^N\setminus\ol{\mc{W}}\,.
\end{cases}
\Ff
\end{definition}

The \ass\ provides one with a description of the upper level sets $E_\lambda(t)$, $\lambda\in(0,1)$, 
for large $t$. Namely, if the initial datum 
$u_0$ is compactly supported then 
\[
\frac{1}{t} E_\lambda(t) \to \W\tinfty,
\]
in the sense of the Hausdorff distance, see \cite[Proposition~2.7]{GHR1}.
As said before, in the homogeneous framework the \ass\ for invading solutions with 
compactly supported initial data is $\W=B_{c^*}$.
In the heterogeneous setting, the \ass, if it exists, is not expected to be symmetric.

The \ass\ encodes the \emph{spreading speed} in any direction 
$e\in\Sph$, which is defined as a quantity $w(e)>0$ satisfying
\Fi{we}
\lim_{t\to+\infty}u(t, cte)=
\begin{cases}
	\displaystyle1 & \text{if }0\leq c<w(e),\vspace{3pt}\\
	\displaystyle0 & \text{if }c>w(e)\,.
\end{cases}
\Ff
One has that
\[
\W=\{ re : e\in\Sph,\ 0\le r < w(e) \}.
\]


\subsection{The Freidlin--G\"artner formula}

Using large deviation probabilistic techniques, Freidlin and G\"artner~\cite{FreidlinGartner} proved the existence 
of the spreading speeds for compactly supported initial data, 
for the Fisher-KPP equation in periodic media. They obtained 
the following formula:
$$
w(e)=\min_{\su{z\in\R^N}{z\.e>0}}\frac{k(z)}{z\.e},
$$
where $k(z)$ is the periodic \pe\ of the linear operator
$$L_z:=\nabla\.(A\nabla)-2z\.A\nabla+q\.z+
\big(-\nabla\.(Az)-q\.z+z\.Az+\partial_u f(x,0)\big).$$
Some of the ideas of \cite{FreidlinGartner} laid the ground for the singular perturbation theory of viscosity solutions,
see e.g.~\cite{ES}.
Several years later, in \cite{BHNperiodic} (see also \cite{W02}), it has been shown 
that, for given $e\in\Sph$, the quantity
$c^*(e):= \min_{\lambda>0}k(\lambda e)/\lambda$
coincides with the critical (or minimal) speed of {\em \PTF s} in the
direction~$e$, hence
Freidlin-G\"artner's formula can be rewritten as
\Fi{FG}
w(e)=\min_{\substack{\xi\in\Sph\\ \xi\.e>0}}\frac{c^*(\xi)}{\xi\.e}.
\Ff
Pulsating traveling fronts,
introduced in parallel ways by Shigesada, Kawasaki and Teramoto
\cite{SKT} and Xin \cite{Xin91-4},
 are the natural generalization of planar traveling fronts to the periodic setting.

\begin{definition}
	A \emph{pulsating traveling front} for \eqref{eq:main}
	in a direction $e\in\Sph$ with a speed $c\in\R$ is a solution of 
	the form
	\[
	u(t,x)=\phi(x, x\cdot e - ct),
	\]
	where $\phi(x,z)$ is periodic in $x$ and satisfies
	\[\phi(x,z)\to 1 \text{ as $z\to -\infty$ \quad and \quad}
	\phi(x,z)\to 0 \text{ as $z\to +\infty$}.\]
\end{definition}

Unlike classical traveling waves in homogeneous media, these fronts oscillate in space due to periodic coefficients. They can be
seen as perturbations of planar fronts, indeed they are {\em almost planar}, 
meaning that their level sets are trapper between
parallel hyperplanes.

The existence of \PTF s in the Monostable, Ignition and --~under suitable assumptions~--
Bistable cases, has been derived in \cite{BerestyckiHamel,Ducrot,GRterrace,W02,Xin91-4,Xin93}.
In the Monostable case, they exist if and only if $c\geq c^*(e)$, whereas in the 
Ignition and Bistable cases they exist if and only if $c=c^*(e)$, for some value $c^*(e)>0$.

In what follows, we refer to $c^*(e)$ as the {\em critical} speed of \PTF s in the direction
$e$, which is the minimal one in the Monostable case, and the unique one in the Ignition-Bistable cases.

One eventually derives from the
Freidlin-G\"artner formula the following expression for the \ass: 
\[
\W=\bigcap_{e\in\Sph} \{x\in\R^N\ :\ x\cdot e < c^*(e)\}\,.
\]
that is, $\W$ is the {\em Wulff shape} of the envelop set of the speeds $c^*(e)$ of pulsating fronts.
One deduces in particular that $\W$ is convex.


\subsection{The \ass\ for general reaction terms}\label{sec:ass}

In the work \cite{RossiFG} it is proved that the
Freidlin-G\"artner formula holds whenever \PTF s exist.
This result recovers the existence of the \ass\ derived by
Weinberger \cite{W02} through the analysis of the Poincar\'e map from the point
of view of discrete dynamical systems. In \cite{W02}, the formula involving
the speeds of fronts is only obtained  in the Monostable case.

We define the Monostable, Ignition, Bistable nonlinearities in the heterogeneous framework as follows:

\begin{align}
\text{\em Monostable }\quad & \forall s\in(0,1),\quad 
\min_{x\in\R^N}f(x,s)\geq0,\quad
\max_{x\in\R^N}f(x,s)>0, \label{mono}\\ \notag \\
\text{\em Ignition } \quad &
\begin{cases}
\exists \theta\in(0,1),\quad\forall (x,s)\in\R^N\times[0,\theta],\quad
f(x,s)=0,\\
\displaystyle\forall s\in(\theta,1),\quad \min_{x\in\R^N}f(x,s)\geq0,\quad
\max_{x\in\R^N}f(x,s)>0,\\
\end{cases}\label{igni}\\ \notag \\
\text{\em Bistable }\quad &
\begin{cases}\exists\sigma\in(0,1),\ \partial_uf(x,\cdot)\leq0
\text{ and $\min_{\R^N} f(\cdot,s)< 0$  in $(0, \sigma]$},\\ 
\text{$\forall e\in\mathbb{S}^{N-1}$, there is a pulsating front with speed $c^*(e)>0$.}
\end{cases}
\end{align}
Observe that the definition of the Bistable case is more implicit, since it assumes the
existence of the pulsating fronts, which allows one to apply Theorems~1.4 and~1.5
of~\cite{RossiFG}. We refer to~\cite{DHZ17,FZ,Ducrot,GRterrace} and the reference therein 
for some sufficient conditions for the existence of pulsating fronts in the Bistable case.

%
%
%
%

\begin{theorem}[\cite{RossiFG}, The generalized Freidlin--G\"artner formula]\label{thm:FG}
	Let the Monostable, Ignition or Bistable condition hold. 
	Let $u$ be an invading solution of \eqref{eq:main}
	with a compactly supported initial datum $u_0$.
	Then $u$ admits an \ass, given~by
	\Fi{W}
	\W=\bigcap_{e\in\Sph} \{x\in\R^N\ :\ x\cdot e < c^*(e)\}\,,
	\Ff
	where $c^*(e)$ is the critical speed of \PTF s in direction~$e$.
%
%
%
\end{theorem}

We point out that, in the KPP case, any nontrivial solution is ``invading'',
while in the other cases, invasion occurs provided the the initial datum $u_0$ is ``not-too-small'', 
see~\cite{Holes,GHR1} for the precise result.

Theorem \ref{thm:FG} extends to general reaction terms the Freidlin-G\"artner's formula~\eqref{FG}
for the spreading speed, as well as the asymptotic description of the upper level sets:
$E_\lambda(t)\sim t\W$ up to a $o(t)$ error.


\subsection{Sketch of the proof of Theorem \ref{thm:FG}}\label{sec:sketch}

We now outline the strategy of the proof of Theorem \ref{thm:FG}.
Define the set $\W$ by \eqref{W}. The goal is to show that $\W$
fulfills the two limits \eqref{ass-cpt} in Definition \ref{ass-cpt}. The first one is a lower estimate of 
the \ass, the second one is an upper estimate. We start with the latter,
which is much simpler.

\subsubsection*{The upper estimate}

Since the constant $1$ is a solution to \eqref{eq:main}, it follows from the 
parabolic strong maximum principle that 
$u(t,x)<1$ for all $t>0$, $x\in\R^N$. Moreover, for any $t>0$,  $u$ has a Gaussian decay in $x$,
because its initial datum is compactly supported. On the other hand, it is known that 
\PTF s decay exponentially.
Then, for any $e\in\Sph$, letting $\phi_e$ be the profile of the front associated with the 
critical speed $c^*(e)$, one can find a quantity $L_e\gg1$
such that
\[\forall t\geq 1,\ x\in\R^N,\qquad
u(t,tx)\leq 
\phi_e(x, x\cdot e - c^*(e)t - L_e)
\]
 (recall that $\phi_e(\.,-\infty)=1$).
One then infers by comparison that, for any $\e>0$,
\[
\begin{split}
\limsup_{t\to+\infty}\Big(\max_{\{x\.e\geq c^*(e)+\e\}}u(t,tx)\Big) &\leq 
\lim_{t\to+\infty}\Big(\max_{\{x\.e\geq c^*(e)+\e\}}
\phi_e(tx, tx\cdot e - c^*(e)t - L_e)\Big)\\
&=\max_{x\in\R^N}\phi_e(x,+\infty)=0.
\end{split}
\]
Since, by a compactness argument, any closed set $C\subset\R^N\setminus\ol{\mc{W}}$
is contained in the union of a finite number of sets of the type $\{x\.e\geq c^*(e)+\e\}$, $e\in\Sph$, $\e>0$, 
one eventually deduces the second limit in \eqref{ass-cpt}.

\subsubsection*{The lower estimate}

The novelty of the approach in \cite{RossiFG} lies in proving the convergence to~$1$ 
simultaneously in all directions, whereas related works in the literature deal with one direction at a time.
This enables one to exploit the geometrical properties of the evolving sets $t\W$.
Indeed, it turns out that these sets evolve with a normal velocity of $c^*(\nu)$, 
where $\nu$ is their outer normal, at least when $\W$ is smooth.
Let us provide some details.

\begin{figure}[H]
\begin{center}
\includegraphics[height=8cm]{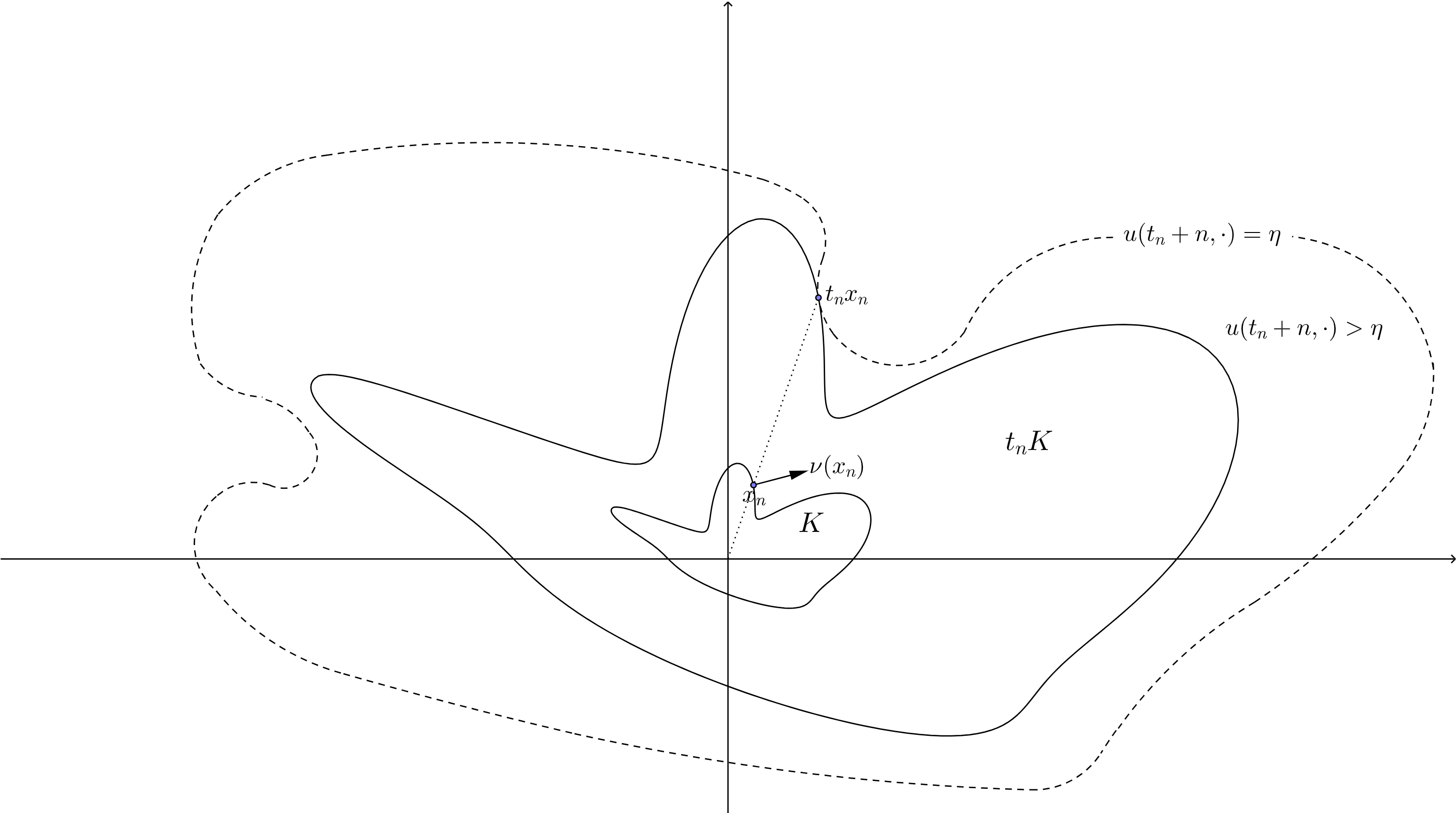}
\caption{The first contact point $(t_n,x_n)$ between $t K$ and the level set 
$\{u(t+n,\.)=\eta\}$.}
\label{fig:K}
\end{center}
\end{figure}

Assume by way of contradiction that 
\eqref{ass-cpt} fails for some compact subset $C$ of~$\mc{W}$,~i.e.
\[\lambda:=\liminf_{t\to+\infty}\,\Big(\min_{x\in C}u(t, tx)\Big)<1.\]
Let us consider another compact set $K$, which is in addition smooth, star-shaped with respect to the origin,
and satisfies
\[C\subset \inter(K)\Subset\W.\]
It follows from the inclusion $C\subset \inter(K)$ 
that
\[\forall n\in\N,\qquad
\liminf_{t\to+\infty}\,\Big(\min_{x\in K}u(t+n, tx)\Big)\leq\lambda.\]
Take an arbitrary $\eta\in(\lambda,1)$. We now look for the first contact time/point between the function
$u(t+n,x)$ and the expanding plateau $\eta\1_{tK}$.
Namely, for any $n\in\N$, there is a time $t_n\geq0$ such that
\[t_n:=\min\{t\geq0\ :\ \min_{x\in K}u(t+n, tx)\leq\eta\}\]
(the above set is nonempty by the contradictory assumption).
In other words, the following hold:
\Fi{tnxn}
\begin{cases}
\forall t\in[0,t_n), \ \forall x\in K,
\quad u(t+n, tx)>\eta\\
\exists x_n\in K,
\quad u(t_n+n, t_nx_n)\leq\eta.
\end{cases}
\Ff
Since $u$ is invading, one infers that $t_n\to+\infty$ as $n\to+\infty$.
Moreover, choosing ${\eta>S}$ from condition \eqref{fS}, one infers from the 
parabolic strong maximum principle that necessarily $x_n\in\partial K$ (otherwise
$u(t_n+n, \.)$ would have a local minimum at $t_nx_n$). 
We denote by $\nu(x_n)$ the outer normal to $K$ at $x_n$.
The situation is depicted in Figure \ref{fig:K}.

We call 
\[u_n(t,x):=u(t_n+n+t, t_nx_n+x).\]
Conditions \eqref{tnxn} rewrite as
\Fi{ntnxn}
\begin{cases}
	\forall t<0, \ \forall x\in (t_n+t)K-\{t_nx_n\},
	\quad u_n(t,x)>\eta\\
	u_n(0,0)=\eta.
\end{cases}
\Ff
Now, the crucial observation is that, being $K\Subset\W$ star-shaped with respect to the origin, the definition \eqref{W} of $\W$ implies that
the family $tK$ evolves by normal velocity less than $c^*(\nu)$, where 
$\nu$ is the outer normal. In addition, since $K$ is (uniformly) smooth, the set $tK$
approaches a half space as $t\to+\infty$, locally around its boundary points, because 
letting~$\rho$ be the radius of the uniform ball conditions
satisfied by $K$, then $tK$ fulfills the ball conditions with the radius
$t\rho$, which tends to infinity as $t\to+\infty$.
Using these properties, one deduces that, for any $z\in\partial K$, there holds
\[
\forall t\in\R,\qquad 
(t_n+t)K - \{t_nz\} \longrightarrow \{x\cdot \nu(z) \leq \sigma(z)t\}\quad\text{ as }\;n\to\infty,
\]
locally uniformly, with $\sigma(z)<c^*(\nu(z))$.
A rigorous derivation of this limit is given in the proof of \cite[Theorem 2.3]{RossiFG}.

As a consequence, it follows from \eqref{ntnxn} that 
the family $\seq{u}$
converges locally uniformly (up to subsequences, by parabolic estimates) toward a function 
$\t u(t,x)$, defined for all $t\in\R$, which satisfies
\Fi{past}
\forall t\leq0, \ \forall x \ :\ x\.\nu\leq \sigma t,
\quad \t u(t,x)\geq\eta,\qquad
\t u(0,0)=\eta,
\Ff
with $\sigma<c^*(\nu)$, where $\nu$ is the limit of (a subsequence of) the exterior normal vectors~$\nu(x_n)$ to $K$.
Moreover, $\t u$ is an entire solution of the same equation~\eqref{eq:main}
(up to a translation). 
But by \eqref{past}, 
its upper level set $\t E_\eta(t)$ moves ``too slowly'' in the past in the direction $\nu$ 
compared with the critical front $\phi_\nu$, namely,
\[\forall t\leq 0,\quad
\t E_\eta(t)\supset\{x\in\R^N\ :\ x\.\nu\leq \sigma t\},\]
with $\sigma<c^*(\nu)$.
Because $\eta>S$, one can indeed show by comparison with $\phi_\nu$ that the above property yields
$\t u\equiv1$, contradicting the equality $\t u(0,0)=\eta$~in \eqref{past}.

This completes the proof of Theorem \ref{thm:FG}.


\section{Convergence towards pulsating traveling fronts}\label{sec:profile}

We now present some recent results obtained in the paper \cite{GHR1} in collaboration with
H.~Guo and F.~Hamel, concerning the asymptotic profile of solutions.

Besides the standing hypotheses of Section \ref{sec:FG},
we assume throughout this section that the states $0$ and $1$ are weakly stable, in the following sense:
\[ 
\exists \theta\in(0,1/2),\quad
s\mapsto f(x,s) \text{ is non-increasing in $[0, \delta]$ and decreasing in $[1-\delta, 1]$}.\]
and moreover:
\[\text{there exists a \PTF\ in any direction $e\in\Sph$, 
with speed $c^*(e)\!>\!0$.}\]
These conditions encompass both the Ignition and the Bistable cases of Section \ref{sec:ass}.

Under these hypotheses, the speed~$c^*(e)$ in any direction~$e$ is unique, and the 
pulsating front is unique up to shift in $t$. We denote it by
\[\Phi_e(t,x)=\phi_e(x,x\. e-c^*(e)t).\]


\subsection{The $\Omega$-limit set}

As mentioned in Section \ref{sec:profile-homo}, in the homogeneous setting, 
solutions emerging from compactly supported initial data converge in profile
to planar traveling fronts. 
The analogous result in the periodic setting is a fundamental open problem in the theory of reaction-diffusion equations. 
To address this question, we considered in \cite{GHR1} the {\em $\O$-limit set} of a solution~$u$, defined by
\be\label{defOmegau}
\Omega(u):=\left\{\begin{array}{rcl}\!\!\psi\in L^\infty(\R^N) & \!\!\!\!:\!\!\!\! & u(t_n,x_n+\.)\displaystyle\mathop{\longrightarrow}_{n\to+\infty}\psi\text{ in $L^\infty_{loc}(\R^N)$,}\\
& \!\!\!\! & \text{for some sequences $(t_n)_{n\in\N}$ in $\R^+$ diverging to $+\infty$}\\
& \!\!\!\! & \text{and $(x_n)_{n\in\N}$ in $\R^N$}\end{array}\!\!\right\}.
\ee
This set is composed by all possible asymptotic profiles of $u$. It can also be written as
\be\label{eqOmegau}
\Omega(u)=\bigcup_{e\in\Sph}\Omega_e(u),
\ee
where
$$\Omega_e(u):=\left\{\begin{array}{rcl}\!\!\psi\in L^\infty(\R^N) & \!\!\!\!:\!\!\!\! & u(t_n,x_n+\.)\displaystyle\mathop{\longrightarrow}_{n\to+\infty}\psi\text{ in $L^\infty_{loc}(\R^N)$ as $n\to+\infty$,}\\
& \!\!\!\! & \text{for some sequences $(t_n)_{n\in\N}$ in $\R^+$ diverging to $+\infty$}\\
& \!\!\!\! & \text{and $(x_n)_{n\in\N}$ in $\R^N\setminus\{0\}$ such that $\frac{x_n}{|x_n|}$ $\displaystyle\mathop{\longrightarrow}_{n\to+\infty}e$}\end{array}\!\!\right\}$$
is called the $\O$-limit set of $u$ in the direction $e$.


\subsection{Compactly supported initial data}

We first focus on compactly supported initial data.

\begin{theorem}[\cite{GHR1}, Corollary 1.7]\label{th1}
Let $u$ be an invading solution to~\eqref{eq:main} with a compactly supported initial datum $u_0$. 
Then, for any $e\in\Sph$, one has 
\be\label{minimizer}\Omega_{e}(u)
\supset
\left\{
\begin{array}{rcl}
\Phi_\xi(t,\.+y) & : & (t,y)\in\R\times\R^N,\ \xi\text{ is a minimizer}\\
& & \text{in the expression of $w(e)$ in \eqref{FG}}
\end{array}\right\}
\ \cup\ \{0,1\}.
\ee
Furthermore, if the boundary of the \ass\ $\W$ (given by \eqref{W}) is everywhere differentiable, then
\be\label{differentiable}
\Omega(u)\supset\big\{\Phi_e(t,\cdot+y):(t,y)\in\R\times\R^N,\ e\in\Sph\big\}\ \cup\ \{0,1\}.
\ee
\end{theorem}

Theorem \ref{th1} asserts that the profile of a \PTF\ appears around the level sets of the solution in any given direction $e$, along some sequence of times. This occurs even in the direction where the \ass\ $\W$ exhibits a corner, i.e.~where the Freidlin-G\"artner formula~\eqref{FG} admits multiple minimizers.
Corners may actually exist,  as shown in~\cite[Corollary~6.2]{GR2}, hence the $\O$-limit set in the corresponding directions contains the profiles of several distinct pulsating~fronts. 
If, on the contrary,~$\partial\W$ is everywhere differentiable then all \PTF s, in every direction, appear as limit profiles
of the solutions, along some sequences of times and~points.

It is natural to wonder whether the equality actually holds in Theorem \ref{th1} when~$\partial\W$ 
is everywhere differentiable, i.e.
$$\Omega(u)=\big\{\Phi_e(t,\cdot+y):(t,y)\in\R\times\R^N,\ e\in\Sph\big\}\ \cup\ \{0,1\}.$$
As we have seen, such conclusion holds for homogeneous equations (where $\W$ is a ball), 
but in the periodic setting it has only been derived in dimension $N=1$, for $f$ either of the Bistable type~\cite{DHZ17,Terrace,Xin91-4,Xin93}, or of the Fisher-KPP type~\cite{HNRR2}. 
To our knowledge, Theorem \ref{th1} is the first convergence-to-front results for periodic equations 
in dimensions higher than one.


\subsection{General initial data}

We finally consider initial data that are not compactly supported.
In such a case, the \ass\ is not determined by just the family of speeds of \PTF s, but also by the initial support of the solution.
For instance, if $u_0=\1_U$ with $U$ being a cone with vertex at the origin, 
and $\partial\W_0$ is of class $C^1$, where $\W_0$
is the \ass\ of invading solutions with compactly supported initial data (given by \eqref{W}), 
then $u$ has the \ass\ $\W=U+\W_0$~\cite{GHR2}.
Actually, due to the possible unboundedness of the \ass, one has to relax the 
Definition~\ref{def:W} by restricting 
to compact subsets $C$
also in the second condition in~\eqref{ass-cpt}.

\begin{theorem}[\cite{GHR1}, Theorem 1.5]\label{thm:GHR1}
Let $u$ be a solution to~\eqref{eq:main} admitting an \ass\ $\W$. 
Let ${z}\in\partial\W$ be a point where $\W$ satisfies the interior and exterior ball conditions, and
let  $\nu$ be the outer normal to~$\W$ at~${z}$, and $\hat z:=z/|z|$. Then the following properties hold:
\begin{enumerate}[$(i)$\;]
\item ${z}\.\nu=c^*(\nu)$\,;
\item $\displaystyle\O_{\hat{z}}(u)\supset\big\{\Phi_\nu(t,\cdot+y):(t,y)\in\R\times\R^N\big\}\cup\{0,1\}$.
\end{enumerate}
\end{theorem}

We point out that under the current hypotheses, cf.~Section \ref{sec:profile},
the invasion property holds for large enough initial data, hence in particular whenever it admits a
nonempty \ass\ $\W$. This implies that $\W$ contains a neighborhood of the origin and thus $z\neq0$ in \thm{GHR1}.

It is easily seen that, in the case of compactly supported initial data, where
$\W$ is given by the Freidlin-G\"artner formula \eqref{W}, at any point~${z}=w(\hat z)\hat z\in\partial\W$ where~$\W$ admits 
a tangent plane, the expression~\eqref{FG} for $w(\hat{z})$ 
is uniquely minimized by $\xi=\nu$, where~$\hat z:=z/|z|$ and~$\nu$ is the outer 
normal to~$\W$ at~$z$, that is,
$${z}\.\nu=c^*(\nu).$$
Therefore, we refer to the above identity as the ``regular Freidlin-G\"artner formula''. 
\thm{GHR1} shows, first, that such formula holds for any solution $u$ at every regular point $z$ of~$\partial\W$ without requiring the boundedness of the support of the initial datum, and second, that the $\O$-limit set of $u$ in the direction $\hat z$ contains profiles of pulsating traveling fronts. Regular here means that $\W$ satisfies both the interior and the exterior ball conditions, 
i.e.~there exists an open ball~$B$ (resp.~$B'$) such that $B\subset\W$ (resp.~$B'\subset\R^N\setminus\W$) and ${z}\in\partial B$ (resp.~${z}\in\partial B'$). The simultaneous validity of the two conditions implies in particular that~$\W$ has a tangent plane at~$z$.

A generalization of Theorem \ref{thm:GHR1} is provided by \cite[Corollary 1.6]{GHR1}
which deals with non-regular points of $\partial\W$
by using the notion of generalized normal set --~i.e.~the set of limits of outer normals of regular approximated~points.


\subsection{Outline of the proofs}

The proof of Theorems \ref{th1} and \ref{thm:GHR1} rely on the same basic idea, which
can be viewed as a variant of Berestycki-Nirenberg's sliding method.
It is also similar to the proof of Theorem \ref{thm:FG} outlined in Section \ref{sec:sketch}.
We recall that there we focused on the first contact time/point between
$u(t+n,x)$ and the expanding plateau $\eta\1_{tK}$.
In the proof of Theorem \ref{thm:GHR1}, 
the goal is to replace the plateau with $\Phi_\nu$, 
the \PTF\ in the direction of the outer normal $\nu$. 
However, in order to do that, one has to make the \PTF\ compactly supported, which is achived by:
\begin{enumerate}
\item making its profile $\phi_\nu(x,z)$ radial in the variable $z$;
\item decreasing $\phi_\nu$, which results in a subsolution for values close to the states $0$ and $1$, thanks to their
stability.
\end{enumerate}
One then focuses on the first contact time/point between 
$u(t+n,x)$ and this perturbed \PTF, lets the perturbation go to $0$ and eventually finds a contact
point between a limit of space/time translations of $u$, and the original \PTF\ $\Phi_\nu$.
One eventually concludes by the strong maximum principle that these two functions coincide, 
that is, $\Phi_\nu\in\Omega_{\hat z}(u)$.

The regularity of $\partial\W$ at $z$ is crucial in this argument 
in order to control the location of the contact points and close the argument.

%


\section{Discussion}

The results presented in this note
show that pulsating traveling fronts are not only special solutions, but fundamental building blocks of the long-time dynamics. Even when the solution does not converge globally to a single wave, its asymptotic behavior is governed by these fronts both locally (profile emergence) and globally (spreading sets).

Theorem \ref{thm:FG}
describes the global geometry of invasion: the solution spreads linearly in time, and the shape of the invaded region is determined by the speeds of \PTF s through their Wulff shape $\W$. 
In homogeneous media, $\mathcal{W}$ would be a ball; in the periodic framework it is generally anisotropic.

Section \ref{sec:profile} is concerned with the convergence towards \PTF s.
The results contained there do not claim full convergence of the solution to a single 
front, but rather the \emph{appearance} of front profiles along subsequences. The global
dynamics might in principle be more complicated.



\end{document}